\documentclass[12pt,a4paper]{amsart}
\usepackage{t1enc}
\usepackage[latin1]{inputenc}
\usepackage{amssymb}
\usepackage{amsthm}
\usepackage{amsmath}
\usepackage{latexsym}
\usepackage{array}
\usepackage{mathrsfs}
\usepackage[all]{xy}
\usepackage{graphics}
\usepackage{xcolor}
\title{{AROUND RATIONALITY OF CYCLES}}
\date{}
\author{{ Raphaël Fino}}
\address
{UPMC Sorbonne Universit\'es\\
Institut de Math\'ematiques de Jussieu\\
Paris\\
FRANCE}
\address
{{\it Web page:}
{\tt www.math.jussieu.fr/\~{ }fino}}
\email {fino {\it at} math.jussieu.fr}
\theoremstyle{definition}
\newtheorem{defi}{Définition}[section]
\newtheorem{rem}[defi]{Remark}
\newtheorem{lemme}[defi]{Lemma}

\newtheorem{prop}[defi]{Proposition}
\newtheorem{thm}[defi]{Theorem}
\begin{document}

\begin{abstract}
In this article we prove certain results comparing rationality of algebraic cycles over the function field of a quadric and over the base field.
Those results have already been proved by Alexander Vishik in the case of characteristic $0$, which allowed him to work with algebraic cobordism theory. Our proofs use the modulo $2$ Steenrod operations in the Chow theory and work in any characteristic $\neq 2$. 
\end{abstract} 

\maketitle

\medskip
  
\medskip

In characteristic $0$, the results of this note (Theorem 1.1 and 2.4 and Proposition 2.1) have been proved several years ago by Alexander Vishik
in [3] (exact references are given right before each statement) with the help of the algebraic cobordism theory and especially \textit{symmetric operations} of [4]. In fact, putting aside characteristic, the original versions are stronger. Indeed, an exponent $2$ element appears in our conclusions, while the use of symmetric operations in the algebraic cobordism theory allows to obtain results without an exponent $2$ element (see Remark on page 370 of [3]). 

In a way, most results of this note are generalizations of some results proved by Nikita Karpenko in [2].
Our proofs are, to a great extent, inspired by the proofs of [2].

In our proofs, the base field is allowed to be of any characteristic different from $2$ because the Landweber-Novikov operations used in [3, Remark after proof of Theorem 3.1] are replaced here by the Steenrod operations on the modulo $2$ Chow groups. 

We refer to [3] and [2] for an introduction into the subject. Notation is introduced in the beginning of Section 1.

\section{Main Result}

Let $F$ be a field of characteristic $\neq 2$, $Q$ a smooth projective quadric over $F$ of dimension $n\geq 0$, $Y$ a smooth quasi-projective 
$F$-variety (a \textit{variety} is a separated scheme of finite type over a field).

The function field  $F(Q)$ is defined if $n\geq 1$ or if $Q$ is anisotropic. In the case of $n=0$ and isotropic $Q$ we have 
$Q= Spec F\coprod Spec F$ and we set $F(Q):=F$.

We write $CH(Y)$ for the integral Chow group of $Y$ (see [1, Chapter X]) and we write $Ch(Y)$ for $CH(Y)$ modulo $2$.
We write $\overline{Y}:=Y_{\overline{F}}$ where $\overline{F}$ is an algebraic closure of $F$.
Let $X$ be a geometrically integral variety over $F$. An element $\overline{y}$ of $Ch(\overline{Y})$ (or of $CH(\overline{Y})$)
is \textit{$F(X)$-rational} if its image $\overline{y}_{\overline{F}(X)}$ under $Ch(\overline{Y})\rightarrow Ch(Y_{\overline{F}(X)})$
(resp. $CH(\overline{Y})\rightarrow CH(Y_{\overline{F}(X)})$) is in the image of $Ch(Y_{F(X)})\rightarrow Ch(Y_{\overline{F}(X)})$ 
(resp. $CH(Y_{F(X)})\rightarrow CH(Y_{\overline{F}(X)})$). Finally, an element $\overline{y}$ of $Ch(\overline{Y})$ (or of $CH(\overline{Y})$)
is called \textit{rational} if it is in the image of $Ch(Y)\rightarrow Ch(\overline{Y})$ (resp. $CH(Y)\rightarrow CH(\overline{Y})$).

\medskip

In a way, the following result is a generalization of [3, Theorem 3.1(1)]. Indeed, the use of the Steenrod operations on the modulo $2$ Chow groups allows to obtain a valid result in any characteristic different from $2$. Nevertheless, an exponent $2$ element appears in our conclusion while it is not the case in [3, Theorem 3.1(1)]. In addition, this result is also a generalization of [2, Theorem 2.1] in the sense that it 
allows a larger codimension for the considered cycle.

\begin{thm}
\textit{Assume that $m<n/2+j$. Let $\overline{y}$ be an $F(Q)$-rational element of $Ch^m(\overline{Y})$.
Then $S^{j}(\overline{y})$ is the sum of a rational element and the class modulo $2$ of an integral element of exponent $2$.}
\end{thm}

\begin{proof}
We assume that $m\geq 0$ in the proof. We also assume that $j \leq m$ (otherwise we get $S^{j}(\overline{y})=0$). The element $\overline{y}$ being $F(Q)$-rational, there exists $y \in Ch^m(Y_{F(Q)})$ mapped to $\overline{y}_{\overline{F}(Q)}$ under the homomorphism 
\[Ch^m(Y_{F(Q)})\rightarrow Ch^m(Y_{\overline{F}(Q)}).\]

Let us fix an element $x \in Ch^m(Q\times Y)$ mapped to $y\;\;mod\;2\;$ under the surjection
\[Ch^m(Q\times Y)\twoheadrightarrow Ch^m(Y_{F(Q)}).\] 
Since over $\overline{F}$ the variety $Q$ becomes cellular, the image $\overline{x}\in Ch^m(\overline{Q}\times \overline{Y})$ 
of $x$ decomposes as
\[\overline{x}=h^0\times y^m +\cdot \cdot \cdot
+h^{[\frac{n}{2}]} \times y^{m-[\frac{n}{2}]}+l_{[\frac{n}{2}]}\times z^{m+[\frac{n}{2}]-n}+
\cdot \cdot \cdot +l_{[\frac{n}{2}]-j}\times z^{m+[\frac{n}{2}]-j-n}\]
with some $y^i \in Ch^i(\overline{Y})$ and some $z^i \in Ch^i(\overline{Y})$, where $y^m=\overline{y}$, and where $h^i \in Ch^i(\overline{Q})$ is the $i$th power of the hyperplane section class while $l_i \in Ch_i(\overline{Q})$ is the class of an $i$-dimensional
subspace of $\mathbb{P}(W)$, where $W$ is a maximal totally isotropic subspace associated with the 
quadric $\overline{Q}$ (see [1, \S68]).

\medskip

For every $i=0,...,m$, let $s^i$ be the image in $CH^{m+i}(\overline{Q}\times \overline{Y})$ of an element in $CH^{m+i}(Q\times Y)$ 
representing $S^i(x)\in Ch^{m+i}(Q\times Y)$. We also set $s^i:=0$ for $i>m$ as well as for $i<0$. 

\medskip

The integer $n$ can be uniquely written in the form $n=2^t-1+s$, where $t$ is a non-negative integer and $0 \leq s < 2^t$.
Let us denote $2^t-1$ as $d$. Since $d \leq n$, we can fix a smooth subquadric $P$ of $Q$ of dimension $d$; we write $in$
for the imbedding 
\[(P\hookrightarrow Q)\times id_Y: P\times Y \hookrightarrow Q\times Y.\]

\begin{lemme}
\textit{For any integer $r$, one has
\[S^r pr_{\ast}in^{\ast}x=\sum_{i=0}^r pr_{\ast}(c_i(-T_P)\cdot in^{\ast}S^{r-i}(x))\;\;\;\;\;in\;\;Ch^{r+m-d}(Y)\]
(where $T_P$ is the tangent bundle of $P$, $c_i$ are the Chern classes, et $pr$ is the projection 
$P\times Y \rightarrow Y$).}
\end{lemme}

\begin{proof}
The morphism $pr:P\times Y \rightarrow Y$ is a smooth projective morphism between smooth schemes. Thus, for any integer $r$, 
we have by [1, Proposition 61.10], 
\[S^r\circ pr_{\ast}=\sum_{i=0}^r pr_{\ast}(c_i(-T_{pr})\cdot S^{r-i})\]
where $T_{pr}$ is the relative tangent bundle of $pr$ over $P\times Y$. Furthermore, since $pr$ is the projection $P\times Y \rightarrow Y$, 
one has $T_{pr}=T_P$. Hence, we get
\[S^r pr_{\ast}in^{\ast}x=\sum_{i=0}^r pr_{\ast}(c_i(-T_P)\cdot S^{r-i}(in^{\ast}x)).\]

Finally, since $in:P\times Y \hookrightarrow Q\times Y$ is a morphism between smooth schemes, the Steenrod operations of cohomological type
commute with $in^{\ast}$ (see [1, Theorem 61.9]), we are done.
\end{proof}

\medskip

We apply Lemma 1.2 taking $r=d+j$. Since $pr_{\ast}in^{\ast}x \in Ch^{m-d}(Y)$ and $m-d<d+j$ (indeed, $m-d<n/2+j-d$ by assumption,
and $n/2<2d$ thanks to our choice of $d$), we have $S^{d+j} pr_{\ast}in^{\ast}x=0$. 

Hence, we have by Lemma 1.2,
\[\sum_{i=0}^{d+j} pr_{\ast}(c_i(-T_P)\cdot in^{\ast}S^{d+j-i}(x))=0\;\;\;\;\;\text{in}\;\;Ch^{m+j}(Y).\]

In addition, for any $i=0,...,d$, by [1, Lemma 78.1] we have $c_i(-T_P)={-d-2 \choose i}\cdot h^i$, where $h^i \in Ch^i(P)$
is the $i$th power of the hyperplane section class, and where the binomial coefficient is considered modulo $2$.
Furthermore, for any $i=0,...,d$, the binomial coefficient ${-d-2 \choose i}={d+i+1 \choose i}$ is odd (because $d$ is a power
of $2$ minus 1, cf. [1, Lemma 78.6]). Moreover, for $i>d$, we have $c_i(-T_P)=0$ in $CH^i(P)$ by definition of Chern classes.
Thus, we get
\[\sum_{i=0}^{d} pr_{\ast}(h^i\cdot in^{\ast}S^{d+j-i}(x))=0\;\;\;\;\;\text{in}\;\;Ch^{m+j}(Y).\]
Therefore, the element
\[\sum_{i=0}^{d} pr_{\ast}(h^i\cdot in^{\ast}s^{d+j-i}) \in CH^{m+j}(\overline{Y})\]
is twice a rational element.

\medskip

Furthermore, for any $i=0,...,d$, we have 
\[pr_{\ast}(h^i\cdot in^{\ast}s^{d+j-i})=pr_{\ast}(in_{\ast}(h^i\cdot in^{\ast}s^{d+j-i}))\]
(the first $pr$ is the projection $P\times Y \rightarrow Y$ while the second $pr$ is the projection 
$Q\times Y \rightarrow Y$). Since $in$ is a proper morphism between smooth schemes, we have by [1, Proposition 56.9],
\[in_{\ast}h^i\cdot in^{\ast}s^{d+j-i}=in_{\ast}h^i\cdot s^{d+j-i}=h^{n-d+i}\cdot s^{d+j-i}\]
and we finally get
\[pr_{\ast}(h^i\cdot in^{\ast}s^{d+j-i})=pr_{\ast}(h^{n-d+i}\cdot s^{d+j-i}).\]

\medskip

Hence, we get that the element
\[\sum_{i=0}^{d}pr_{\ast}(h^{n-d+i}\cdot s^{d+j-i})\in CH^{m+j}(\overline{Y})\]
is twice a rational element.

\medskip

We would like to compute the sum obtained modulo $4$. Since $s^{d+j-i}=0$ if $d+j-i>m$, the $i$th summand is $0$ for any 
$i<d+j-m$ ($(j-m) \leq 0$ by assumption). Otherwise -- if $i \geq d+j-m$ -- the factor $h^{n-d+i}$ is divisible by $2$ (indeed, we have
$h^{n-d+i}=2l_{d-i}$ because $n-d+i\geq n+j-m>n/2$) and in order to compute the $i$th summand modulo $4$ it suffices to compute
$s^{d+j-i}$ modulo $2$, that is, to compute $S^{d+j-i}(\overline{x})$.

We recall that 
\[\overline{x}=h^0\times y^m +\cdot \cdot \cdot+
h^{[\frac{n}{2}]} \times y^{m-[\frac{n}{2}]}+l_{[\frac{n}{2}]}\times z^{m+[\frac{n}{2}]-n}+
\cdot \cdot \cdot +l_{[\frac{n}{2}]-j}\times z^{m+[\frac{n}{2}]-j-n}.\]
Therefore, we have
\[S^{d+j-i}(\overline{x})=\sum_{k=0}^{[\frac{n}{2}]} S^{d+j-i}(h^k \times y^{m-k})+\sum_{k=0}^{j} S^{d+j-i}(l_{[\frac{n}{2}]-k} 
\times z^{m+[\frac{n}{2}]-k-n}).\]
And we set
\[A_i:=\sum_{k=0}^{[\frac{n}{2}]} S^{d+j-i}(h^k \times y^{m-k})\;\;and\;\;B_i:=\sum_{k=0}^{j} S^{d+j-i}(l_{[\frac{n}{2}]-k} 
\times z^{m+[\frac{n}{2}]-k-n}).\]

\medskip

For any $k=0,...,[\frac{n}{2}]$, we have by [1, Theorem 61.14],
\[S^{d+j-i}(h^k \times y^{m-k})=\sum_{l=0}^{d+j-i}S^{d+j-i-l}(h^k)\times S^l(y^{m-k}).\]
Moreover, for any  $l=0,...,d+j-i$, we have by [1, Corollary 78.5],
\[S^{d+j-i-l}(h^k)={k \choose d+j-i-l}h^{d+j+k-i-l}.\]

Thus, choosing an integral representative $\varepsilon_{k,l}\in CH^{m-k+l}(\overline{Y})$ of $S^l(y^{m-k})$ 
(we choose $\varepsilon_{k,l}=0$ if $l>m-k$), we get that the element 
\[\sum_{k=0}^{[\frac{n}{2}]} \sum_{l=0}^{d+j-i}{k \choose d+j-i-l}(h^{d+j+k-i-l}\times \varepsilon_{k,l}) \in CH^{m+d+j-i}(\overline{Q}\times \overline{Y})\]
is an integral representative of $A_i$.

\medskip

Therefore, for any $i \geq d+j-m$, choosing an integral representative $\tilde{B_i}$ of $B_i$, there exists $\gamma_i \in CH^{m+d+j-i}(\overline{Q}\times \overline{Y})$ such that
\[s^{d+j-i}=\sum_{k=0}^{[\frac{n}{2}]} \sum_{l=0}^{d+j-i}{k \choose d+j-i-l}(h^{d+j+k-i-l}\times \varepsilon_{k,l})
+\tilde{B_i}+ 2\gamma_i.\]
Hence, according to the multiplication rules in the ring $CH(\overline{Q})$ described in [1, Proposition 68.1], 
for any $i \geq d+j-m$, we have 
\[h^{n-d+i}\cdot s^{d+j-i}=2\sum_{k=0}^{[\frac{n}{2}]} \sum_{l=0}^{d+j-i}{k \choose d+j-i-l}(l_{l-j-k}\times \varepsilon_{k,l})
+h^{n-d+i}\cdot \tilde{B_i}+ 4l_{d-i}\cdot \gamma_i. \]

If $k \leq d-i$, one has $j+k\leq d+j-i$, and for any $0\leq l \leq d+j-i$, we have by dimensional reasons, 
\[pr_{\ast}(l_{l-j-k}\times \varepsilon_{k,l})=\left\{\begin{array}{rcl} &\varepsilon_{k,l} & \text{if}\,\,l=j+k\\
                                   &0&\text{otherwise.}
\end{array}\right.\] 
Otherwise $k>d-i$, and $pr_{\ast}(l_{l-j-k}\times \varepsilon_{k,l})=0$ for any $0\leq l \leq d+j-i$. 
Moreover, for $k>d-i$, one has $j+k>j+d-i\geq m >m-k$, therefore $\varepsilon_{k,j+k}=0$.

Thus we deduce the identity
\[pr_{\ast}\left(2\sum_{k=0}^{[\frac{n}{2}]} \sum_{l=0}^{d+j-i}{k \choose d+j-i-l}(l_{l-j-k}\times \varepsilon_{k,l})\right)=
2\sum_{k=0}^{[\frac{n}{2}]}{k \choose d-i-k}\varepsilon_{k,j+k}.\]

Then,
\begin{multline*}
\sum_{i=d+j-m}^d pr_{\ast}\left(2\sum_{k=0}^{[\frac{n}{2}]} \sum_{l=0}^{d+j-i}{k \choose d+j-i-l}(l_{l-j-k}\times \varepsilon_{k,l})\right)\\
  =2\sum_{i=d+j-m}^d \sum_{k=0}^{[\frac{n}{2}]}{k \choose d-i-k}\varepsilon_{k,j+k}.
\end{multline*}

In the latest expression, for every $k=0,...,[\frac{n}{2}]$, the total coefficient near $\varepsilon_{k,j+k}$ is 
\[2\sum_{i=d+j-m}^{d}{k \choose d-i-k}=2\sum_{i=d-2k}^{d-k}{k \choose d-i-k}=2\sum_{s=0}^{k}{k \choose s}=2^{k+1},\]
which is divisible by $4$ for $k \geq 1$.

\medskip

Therefore, the cycle $\sum_{i=d+j-m}^{d}pr_{\ast}(h^{n-d+i}\cdot s^{d+j-i})\in CH^{m+j}(\overline{Y})$ is congruent modulo 4 to
\[2\varepsilon_{0,j}+\sum_{i=d+j-m}^{d}pr_{\ast}(h^{n-d+i}\cdot \tilde{B_i}).\]

Thus, the cycle $2\varepsilon_{0,j}+\sum_{i=d+j-m}^{d}pr_{\ast}(h^{n-d+i}\cdot \tilde{B_i})$ is congruent modulo $4$ to twice a rational element.

\medskip

Finally, the following lemma will lead to the conclusion.

\begin{lemme}
\textit{For any $d+j-m\leq i \leq d$, one can choose an integral representative $\tilde{B_i}$ of $B_i$ so that  
 \[pr_{\ast}(h^{n-d+i}\cdot \tilde{B_i})=0.\]}
\end{lemme}

\begin{proof}
We recall that $B_i:=\sum_{k=0}^{j} S^{d+j-i}(l_{[\frac{n}{2}]-k}\times z^{m+[\frac{n}{2}]-k-n})$.
For any $k=0,...,j$, we have by [1, Theorem 61.14],
\[S^{d+j-i}(l_{[\frac{n}{2}]-k} \times z^{m+[\frac{n}{2}]-k-n})=\sum_{l=0}^{d+j-i}S^{d+j-i-l}(l_{[\frac{n}{2}]-k})\times S^l(z^{m+[\frac{n}{2}]-k-n}).\]
And for any $l=0,...,d+j-i$, we have by [1, Corollary 78.5],
\[S^{d+j-i-l}(l_{[\frac{n}{2}]-k})={n+1-[\frac{n}{2}]+k \choose d+j-i-l}l_{[\frac{n}{2}]-k-d-j+i+l}.\]

Thus, choosing an integral representative $\delta_{k,l}\in CH^{m-k+l}(\overline{Y})$ of $S^l(z^{m+[\frac{n}{2}]-k-n})$ 
(we choose $\delta_{k,l}=0$ if $l>m+[\frac{n}{2}]-k-n$), we get that the element
\[\sum_{k=0}^{j} \sum_{l=0}^{d+j-i}{n+1-[\frac{n}{2}]+k \choose d+j-i-l}(l_{[\frac{n}{2}]-k-d-j+i+l}\times \delta_{k,l}) \in CH^{m+d+j-i}(\overline{Q}\times \overline{Y})\]
is an integral representative of $B_i$. Let us note it $\tilde{B_i}$.

\medskip

Hence, we have
\[h^{n-d+i}\cdot \tilde{B_i}=\sum_{k=0}^{j} \sum_{l=0}^{d+j-i}{n+1-[\frac{n}{2}]+k \choose d+j-i-l}(l_{[\frac{n}{2}]-k-n-j+l}\times \delta_{k,l}).\]

Moreover, we have 
\[pr_{\ast}(l_{[\frac{n}{2}]-k-n-j+l}\times \delta_{k,l})\neq 0\Longrightarrow l=j+k+n-[\frac{n}{2}].\]

Furthermore, for any $0\leq k \leq j$, we have $d+j-i \leq m <j+\frac{n}{2}\leq j+n-[\frac{n}{2}]\leq j+k+n-[\frac{n}{2}]$.
Thus, for any $0 \leq l\leq d+j-i$ and for any $0\leq k \leq j$, we have $pr_{\ast}(l_{[\frac{n}{2}]-k-n-j+l}\times \delta_{k,l})=0$.
It follows that $pr_{\ast}(h^{n-d+i}\cdot \tilde{B_i})=0$ and we are done. 
\end{proof}

We deduce from Lemma 1.3 that the cycle $2\varepsilon_{0,j} \in CH^{m+j}(\overline{Y})$ is congruent modulo $4$ to twice a rational cycle.
Therefore, there exist a cycle $\gamma \in CH^{m+j}(\overline{Y})$ and a rational cycle $\alpha \in CH^{m+j}(\overline{Y})$ so that
\[2\varepsilon_{0,j}=2\alpha + 4\gamma,\]
hence, there exists an exponent $2$ element $\delta \in CH^{m+j}(\overline{Y})$ so that
\[\varepsilon_{0,j}=\alpha + 2\gamma + \delta.\]

\medskip

Finally, since $\varepsilon_{0,j}$ is an integral representative of $S^{j}(\overline{y})$, we get that
$S^{j}(\overline{y})$ is the sum of a rational element and the class modulo $2$ of an integral element of exponent $2$. 
We are done with the proof of Theorem 1.1.
\end{proof}

\section{Other results}

In this section we continue to use notation introduced in the beginning of Section~1.
In the same way as before, the following proposition is a generalization of [3, Proposition 3.3(2)] 
(although, putting aside characteristic, our proposition is still weaker than the original version 
in the sense that an exponent $2$ element appears in the conclusion).

\begin{prop}
\textit{Let $x \in Ch^m(Q\times Y)$ be some element, and $y^i, z^i \in Ch^i(\overline{Y})$ be the coordinates of $\overline{x}$ as 
in the beginning of proof of Theorem 1.1.
Assume that $m=[\frac{n+1}{2}]+j$ and that $n\geq 1$.} 
\textit{Then $S^j(y^m)+y^m\cdot z^j$ differs from a rational element by the class of an exponent $2$ element of 
$CH^{m+j}(\overline{Y})$.}
\end{prop}

\begin{proof}
The image $\overline{x}\in Ch^m(\overline{Q}\times \overline{Y})$ of $x$ decomposes as
\[\overline{x}=h^0\times y^m +\cdot \cdot \cdot
+h^{[\frac{n}{2}]} \times y^{m-[\frac{n}{2}]}+l_{[\frac{n}{2}]}\times z^{m+[\frac{n}{2}]-n}+
\cdot \cdot \cdot +l_{[\frac{n}{2}]-j-1}\times z^{m+[\frac{n}{2}]-j-n}.\] 

\medskip

Let $\textbf{x} \in CH^m(Q \times Y)$ be an integral representative of $x$. 
The image $\overline{\textbf{x}} \in CH^m(\overline{Q} \times \overline{Y})$ decomposes as 
\[\overline{\textbf{x}}=h^0\times \textbf{y}^m +\cdot \cdot \cdot
+h^{[\frac{n}{2}]} \times \textbf{y}^{m-[\frac{n}{2}]}+l_{[\frac{n}{2}]}\times \textbf{z}^{m+[\frac{n}{2}]-n}+
\cdot \cdot \cdot +l_{[\frac{n}{2}]-j-1}\times \textbf{z}^{m+[\frac{n}{2}]-j-n}\]
where the elements $\textbf{y}^i \in CH^i(\overline{Y})$ (resp. $\textbf{z}^i \in CH^i(\overline{Y})$) are some integral representatives of 
the elements $y^i$ (resp. $z^i$).

\medskip

For every $i=0,...,m-1$, let $s^i$ be the image in $CH^{m+i}(\overline{Q}\times \overline{Y})$ of an element in $CH^{m+i}(Q\times Y)$ 
representing $S^i(x)\in Ch^{m+i}(Q\times Y)$. We also set $s^i:=0$ for $i>m$ as well as for $i<0$. Finally, we set $s^0:=\overline{\textbf{x}}$ 
and $s^m:={(s^0)}^2$. Therefore, for any integer $i$, $s^i$ is the image in $CH^{m+i}(\overline{Q} \times \overline{Y})$ of an integral representative of $S^i(x)$.

\medskip

The integer $n$ can be uniquely written in the form $n=2^t-1+s$, where $t$ is a non-negative integer and $0 \leq s < 2^t$.
Let us denote $2^t-1$ as $d$.

\medskip

We would like to use again Lemma 1.2 to get that the sum
\[\sum_{i=d+j-m}^{d}pr_{\ast}(h^{n-d+i}\cdot s^{d+j-i})\in CH^{m+j}(\overline{Y})\]
is twice a rational element. To do this, it suffices to check that $m-d<d+j$. 
Then the same reasoning as the one used during the proof of Theorem 1.1 gives us the desired result.

We have $m-d=[\frac{n+1}{2}]+j-d=d+j+([\frac{n+1}{2}]-2d)$, and since our choice of $d$ and the assumption $n\geq 1$,
one can easily check that $2d>[\frac{n+1}{2}]$. Thus we do get that the sum 
\[\sum_{i=d+j-m}^{d}pr_{\ast}(h^{n-d+i}\cdot s^{d+j-i})\in CH^{m+j}(\overline{Y})\]
is twice a rational element. We would like to compute that sum modulo $4$.

\medskip

For any $i\geq d+j-m$, the factor $s^{d+j-i}$ present in the $i$th summand is congruent modulo $2$ to 
$S^{d+j-i}(\overline{x})$, which is represented by $\tilde{A_i}+ \tilde{B_i}$, where
\[\tilde{A_i}:=\sum_{k=0}^{[\frac{n}{2}]} \sum_{l=0}^{d+j-i}{k \choose d+j-i-l}(h^{d+j+k-i-l}\times \varepsilon_{k,l})\]
and
\[\tilde{B_i}:=\sum_{k=0}^{j} \sum_{l=0}^{d+j-i}{n+1-[\frac{n}{2}]+k \choose d+j-i-l}(l_{[\frac{n}{2}]-k-d-j+i+l}\times \delta_{k,l})\]
where $\varepsilon_{k,l}\in CH^{m-k+l}(\overline{Y})$ (resp. $\delta_{k,l}\in CH^{m-k+l}(\overline{Y})$) is an integral representative of $S^l(y^{m-k})$ (resp. of $S^l(z^{m+[\frac{n}{2}]-k-n})$), 
and we choose $\varepsilon_{k,l}=0$ if $l>m-k$ (resp. $\delta_{k,l}=0$ if $l>m+[\frac{n}{2}]-k-n$).
Finally, in the case of even $m-j$ , we choose $\varepsilon_{\frac{m-j}{2},\frac{m+j}{2}}=(\textbf{y}^{\frac{m+j}{2}})^2$.

\medskip

Furthermore, for any $i\geq d+j-m$, we have 
\[h^{n-d+i}\cdot \tilde{B_i}=\sum_{k=0}^{j} \sum_{l=0}^{d+j-i}{n+1-[\frac{n}{2}]+k \choose d+j-i-l}(l_{[\frac{n}{2}]-k-n-j+l}\times \delta_{k,l}).\] 

And we have 
\[pr_{\ast}(l_{[\frac{n}{2}]-k-n-j+l}\times \delta_{k,l})\neq 0\Longrightarrow l=j+k+n-[\frac{n}{2}].\]

On the one hand, for any $i>d+j-m$, we have $d+j-i<m=n-[\frac{n}{2}]+j \leq j+k+n-[\frac{n}{2}]$.
Hence, for any $0 \leq l\leq d+j-i$ and for any $0\leq k \leq j$, we have $pr_{\ast}(l_{[\frac{n}{2}]-k-n-j+l}\times \delta_{k,l})=0$.
Then, for any $i>d+j-m$, we get that $pr_{\ast}(h^{n-d+i}\cdot \tilde{B_i})=0$. 

\medskip

On the other hand, for $i=d+j-m$, we have $d+j-i=j+n-[n/2]$ and 
\[l=j+k+n-[\frac{n}{2}] \Longleftrightarrow k=0\;\;and\;\;l=d+j-i.\]
Thus, we have 
\[pr_{\ast}(h^{n+j-m}\cdot \tilde{B}_{d+j-m})=\delta_{0,m}.\]
Since $m>m+[n/2]-n$, we get that $\delta_{0,m}=0$.

\medskip

Therefore, for any $i\geq d+j-m$, we have  
\[pr_{\ast}(h^{n-d+i}\cdot \tilde{B_i})=0.\]

\medskip

Then, for any $i>d+j-m$, the cycle $h^{n-d+i}$ is divisible by $2$. Hence, according to the multiplication rules in the ring $CH(\overline{Q})$ described in [1, Proposition 68.1] and by doing the same computations as those done during the proof of Theorem 1.1, for any $i>d+j-m$, we get the congruence
\[pr_{\ast}(h^{n-d+i}\cdot s^{d+j-i})\equiv 2\sum_{k=0}^{[\frac{n}{2}]}{k \choose d-i-k}\varepsilon_{k,j+k} \;\;(mod\;4).\]
Moreover, since $d-i-k\leq k$ if and only if $k\leq [\frac{m-j}{2}]$, for any $i>d+j-m$, we have the congruence
\begin{equation} pr_{\ast}(h^{n-d+i}\cdot s^{d+j-i})\equiv 2\sum_{k=0}^{[\frac{m-j}{2}]}{k \choose d-i-k}\varepsilon_{k,j+k} \;\;(mod\;4). \end{equation}

Now, we would like to study the $(d+j-m)$th summand, that is to say the cycle $pr_{\ast}(h^{n+j-m}\cdot s^m)$ modulo $4$.
That is the purpose of the following lemma.

\medskip

\begin{lemme}
\textit{One has
\[pr_{\ast}(h^{n+j-m}\cdot s^m)\equiv \left\{\begin{array}{rcl} &2\varepsilon_{\frac{m-j}{2},\frac{m+j}{2}}+2\textbf{y}^m\cdot \textbf{z}^j 
 \;\;(mod\;4) & if\;m-j\;is\;even\\
                                   &2\textbf{y}^m\cdot \textbf{z}^j  \;\;(mod\;4) &if\;m-j\;is\;odd.
\end{array}\right.\]} 
\end{lemme}

\begin{proof}
We recall that $s^m=(\overline{\textbf{x}})^2$. Thus, we have
\[h^{n+j-m}\cdot s^m=h^{n+j-m}\cdot (A+B+C)\]
where
\[A:=\sum_{0\leq i,l \leq [\frac{n}{2}]} h^{i+l}\times (\textbf{y}^{m-i}\cdot \textbf{y}^{m-l}),\]
\[B:=\sum_{0\leq i,l \leq j} (l_{[\frac{n}{2}]-i}\cdot l_{[\frac{n}{2}]-l})\times (\textbf{z}^{j-i}\cdot \textbf{z}^{j-l})\]
and
\[C:=2\sum_{i=0}^{[\frac{n}{2}]} h^i\times \textbf{y}^{m-i} \cdot \sum_{l=0}^j l_{[\frac{n}{2}]-l}\times \textbf{z}^{j-l}.\]

First of all, we have
\[h^{n+j-m}\cdot A=\sum_{0\leq i,l \leq [\frac{n}{2}]} h^{n+j-m+i+l}\times (\textbf{y}^{m-i}\cdot \textbf{y}^{m-l}).\]
Now we have $m=[\frac{n+1}{2}]+j$, so $n+j-m+i+l=[\frac{n}{2}]+i+l$. Thus, if $i\geq 1$ or $l\geq 1$, we have 
$n+j-m+i+l>[\frac{n}{2}]$, and in this case we have $h^{n+j-m+i+l}=2l_{m-i-l-j}$. Therefore, 
the cycle $h^{n+j-m}\cdot A$ is equal to
\[h^{n+j-m}\times (\textbf{y}^m)^2 + 4\sum_{1\leq i,l \leq [\frac{n}{2}]}l_{m-i-l-j} \times (\textbf{y}^{m-i}\cdot \textbf{y}^{m-l})
+2\sum_{i=1}^{[\frac{n}{2}]} l_{m-j-2i}\times (\textbf{y}^{m-i})^2.\]

\medskip

Then, since $n\geq 1$, we have $n+j-m \neq n$. It follows that $pr_{\ast}(h^{n+j-m}\times (\textbf{y}^m)^2)=0$.

\medskip

Furthermore, we have 
\[pr_{\ast}(\sum_{i=1}^{[\frac{n}{2}]} l_{m-j-2i}\times (\textbf{y}^{m-i})^2)=\left\{\begin{array}{rcl} &(\textbf{y}^{\frac{m+j}{2}})^2 & \text{if $m-j$ is even}\\
                                   &0& \text{if $m-j$ is odd.}
\end{array}\right.\] 
Therefore, $pr_{\ast}(h^{n+j-m}\cdot A)$ is congruent modulo $4$ to $2\varepsilon_{\frac{m-j}{2},\frac{m+j}{2}}$ if $m-j$ is even, and to $0$ if $m-j$ is odd.

\medskip

Then, by dimensional reasons, we have $l_{[\frac{n}{2}]-i}\cdot l_{[\frac{n}{2}]-l}=0$ if $i\geq 1$ or if $l\geq 1$.
Hence, we have $B=(l_{[\frac{n}{2}]}\cdot l_{[\frac{n}{2}]})\times (\textbf{z}^j)^2$. It follows that 
\[h^{n+j-m}\cdot B=(l_0 \cdot l_{[\frac{n}{2}]})\times (\textbf{z}^j)^2\]
and $l_0 \cdot l_{[\frac{n}{2}]}=0$ by dimensional reasons. Therefore, we get that $h^{n+j-m}\cdot B=0$.

\medskip

Finally, we have 
\[h^{n+j-m}\cdot C=2\sum_{i=0}^{[\frac{n}{2}]} h^{n+j-m+i}\times \textbf{y}^{m-i} \cdot \sum_{l=0}^j l_{[\frac{n}{2}]-l}\times \textbf{z}^{j-l}.\]
Now for any $i\geq 1$, we have $n+j-m+i>[\frac{n}{2}]$, and in this case the cycle $h^{n+j-m+i}$ is divisible by $2$. 
Thus, the element $h^{n+j-m}\cdot C$ is congruent modulo $4$ to 
\[2\sum_{l=0}^j (h^{[\frac{n}{2}]}\cdot l_{[\frac{n}{2}]-l})\times (\textbf{y}^{m} \cdot \textbf{z}^{j-l}),\]
and, by dimensional reasons, in the latest sum, each summand is $0$ except the one corresponding to $l=0$.
Therefore, the cycle $h^{n+j-m}\cdot C$ is congruent modulo $4$ to $2l_0 \times (\textbf{y}^{m} \cdot \textbf{z}^{j})$.
It follows that $pr_{\ast}(h^{n+j-m}\cdot C)$ is congruent modulo $4$ to $2\textbf{y}^{m} \cdot \textbf{z}^{j}$.
We are done.
\end{proof}

By the congruence (1) and Lemma 2.2, we deduce that the cycle 
\[\sum_{i=d+j-m}^{d}pr_{\ast}(h^{n-d+i}\cdot s^{d+j-i})\]
is congruent modulo $4$ to
\[2\sum_{i=d+j-m}^d \sum_{k=0}^{[\frac{m-j}{2}]}{k \choose d-i-k}\varepsilon_{k,j+k} + 2\textbf{y}^m\cdot \textbf{z}^j.\]
It follows that the cycle 
\[2\sum_{i=d+j-m}^d \sum_{k=0}^{[\frac{m-j}{2}]}{k \choose d-i-k}\varepsilon_{k,j+k} + 2\textbf{y}^m\cdot \textbf{z}^j\]
is congruent modulo $4$ to twice a rational element $\alpha \in CH^{m+j}(\overline{Y})$.
Then, we finish as in the proof of Theorem 1.1. For every $k=0,...,[(m-j)/2]$, the total coefficient near $\epsilon_{k,j+k}$ is $2^{k+1}$, which is divisible by $4$ for $k \geq 1$.
Therefore, there exists a cycle $\gamma \in CH^{m+j}(\overline{Y})$ such that 
\[2\varepsilon_{0,j}+2\textbf{y}^m\cdot \textbf{z}^j=2\alpha + 4\gamma,\]
hence, there exists an exponent $2$ element $\delta \in CH^{m+j}(\overline{Y})$ so that
\[\varepsilon_{0,j}+\textbf{y}^m\cdot \textbf{z}^j=\alpha + 2\gamma + \delta.\]

\medskip

Finally, since $\varepsilon_{0,j}$ is an integral representative of $S^{j}(y^m)$ and $\textbf{y}^m$ (resp. $\textbf{z}^j$)
is an integral representative of $y^m$ (resp. of $z^j$), we get that $S^j(y^m)~+~y^m\cdot ~z^j$ 
differs from a rational element by the class of an exponent $2$ element of 
$CH^{m+j}(\overline{Y})$.
We are done with the proof of Proposition 2.1.
\end{proof}

\medskip

\begin{rem}
In the case of $j=0$, and if we make the extra assumption that the image of $x$ under the composition
\[Ch^m(Q \times Y) \rightarrow Ch^m(Q_{F(Y)}) \rightarrow Ch^m(Q_{\overline{F}(Y)}) \rightarrow Ch^m(\overline{Q})\]
(the last passage is given by the inverse of the change of field isomorphism)
is rational, then we get the stronger result that the cycle $y^m$ differs from a rational element by the class of an exponent $2$ element of 
$CH^{m}(\overline{Y})$. That is the object of [2, Proposition 4.1]. 
\end{rem}

\medskip

Finally, the following theorem is a consequence of Proposition 2.1. In the same way as before, it is a generalization of [3, Theorem 3.1(2)].
For a variety $X$, we write $r$ for the restriction map $Ch^{\ast}(X)\rightarrow Ch^{\ast}(\overline{X})$.

\begin{thm}
\textit{Assume that $m=[\frac{n+1}{2}]+j$ and $n\geq 1$. Let $\overline{y}$ be an $F(Q)$-rational element of $Ch^m(\overline{Y})$.
Then there exists a rational element $z \in Ch^j(\overline{Y})$ such that $S^{j}(\overline{y})+\overline{y}\cdot z$ is the sum of a rational element and the class modulo $2$ of an integral element of exponent $2$.}
\end{thm}

\begin{proof}
The element $\overline{y}$ being $F(Q)$-rational, there exists $x \in Ch^m(Q\times Y)$ mapped to $\overline{y}_{\overline{F}(Q)}$ under the composition
\[Ch^m(Q \times Y) \rightarrow Ch^m(Y_{F(Q)}) \rightarrow Ch^m(Y_{\overline{F}(Q)}).\]

Moreover, the image $\overline{x}\in Ch^m(\overline{Q}\times \overline{Y})$ of $x$ decomposes as
\[\overline{x}=h^0\times y^m +\cdot \cdot \cdot
+h^{[\frac{n}{2}]} \times y^{m-[\frac{n}{2}]}+l_{[\frac{n}{2}]}\times z^{m+[\frac{n}{2}]-n}+
\cdot \cdot \cdot +l_{[\frac{n}{2}]-j-1}\times z^{m+[\frac{n}{2}]-j-n}\]
with some $y^i \in Ch^i(\overline{Y})$, and some $z^i \in Ch^i(\overline{Y})$, and where $y^m=\overline{y}$. 

Thus, by Proposition 2.1, the cycle $S^{j}(\overline{y})+\overline{y}\cdot z^j$
is the sum of a rational element and the class an element of exponent $2$.

\medskip

Finally, we have by [1,  Proposition 49.20],
\[r\circ(pr)_{\ast}(x\cdot h^{[\frac{n}{2}]})=pr_{\ast}(\overline{x}\cdot h^{[\frac{n}{2}]})=z^j.\]

We are done with the proof of Theorem 2.4.
\end{proof}


\begin{thebibliography}{2}
\bibitem[1]{1} R.~Elman, N.~Karpenko, A.~Merkurjev, \emph{The algebraic and geometric theory of quadratic forms},
American Mathematical Society Colloquium Publications, 56. American Mathematical Society, Providence, RI, 2008.
\bibitem[2]{2} N.~Karpenko, Variations on a theme of rationality of cycles. \textit{Linear Algebraic Groups and Related Structures (preprint server) 443} (2011, Sep 21), 13 pages.
\bibitem[3]{3} A.~Vishik, Generic points of quadrics and Chow groups, \textit{Manuscripta Math. 122, 3 (2007), 365-374}.
\bibitem[4]{4} A.~Vishik, Symmetric operations in algebraic cobordism. \textit{Adv. Math. 213}, 2 (2007), 489-552.
\end{thebibliography}
\end{document}